\documentclass[12pt]{amsart}
\usepackage{times}
\usepackage [latin1]{inputenc}
\begin{document}
\numberwithin{equation}{section}

\def\1#1{\overline{#1}}
\def\2#1{\widetilde{#1}}
\def\3#1{\widehat{#1}}
\def\4#1{\mathbb{#1}}
\def\5#1{\frak{#1}}
\def\6#1{{\mathcal{#1}}}

\newcommand{\de}{\partial}
\newcommand{\R}{\mathbb R}
\newcommand{\al}{\alpha}
\newcommand{\tr}{\widetilde{\rho}}
\newcommand{\tz}{\widetilde{\zeta}}
\newcommand{\tv}{\widetilde{\varphi}}
\newcommand{\hv}{\hat{\varphi}}
\newcommand{\tu}{\tilde{u}}
\newcommand{\tF}{\tilde{F}}
\newcommand{\debar}{\overline{\de}}
\newcommand{\Z}{\mathbb Z}
\newcommand{\C}{\mathbb C}
\newcommand{\Po}{\mathbb P}
\newcommand{\zbar}{\overline{z}}
\newcommand{\G}{\mathcal{G}}
\newcommand{\So}{\mathcal{S}}
\newcommand{\Ko}{\mathcal{K}}
\newcommand{\U}{\mathcal{U}}
\newcommand{\B}{\mathbb B}
\newcommand{\oB}{\overline{\mathbb B}}
\newcommand{\Cur}{\mathcal D}
\newcommand{\Dis}{\mathcal Dis}
\newcommand{\Levi}{\mathcal L}
\newcommand{\SP}{\mathcal SP}
\newcommand{\Sp}{\mathcal Q}
\newcommand{\Ma}{\mathcal M}
\newcommand{\Co}{\mathcal C}
\newcommand{\Hol}{{\sf Hol}(\mathbb D, D)}
\newcommand{\Aut}{{\sf Aut}(\mathbb D)}
\newcommand{\D}{\mathbb D}
\newcommand{\oD}{\overline{\mathbb D}}
\newcommand{\oX}{\overline{X}}
\newcommand{\loc}{L^1_{\rm{loc}}}
\newcommand{\la}{\langle}
\newcommand{\ra}{\rangle}
\newcommand{\thh}{\tilde{h}}
\newcommand{\N}{\mathbb N}
\newcommand{\kd}{\kappa_D}
\newcommand{\Hr}{\mathbb H}
\newcommand{\ps}{{\sf Psh}}

\newcommand{\subh}{{\sf subh}}
\newcommand{\harm}{{\sf harm}}
\newcommand{\ph}{{\sf Ph}}
\newcommand{\tl}{\tilde{\lambda}}
\newcommand{\ts}{\tilde{\sigma}}

\def\v{\varphi}
\def\Re{{\sf Re}\,}
\def\Im{{\sf Im}\,}

\def\dist{{\rm dist}}
\def\const{{\rm const}}
\def\rk{{\rm rank\,}}
\def\id{{\sf id}}
\def\aut{{\sf aut}}
\def\Aut{{\sf Aut}}
\def\CR{{\rm CR}}
\def\GL{{\sf GL}}
\def\U{{\sf U}}

\def\la{\langle}
\def\ra{\rangle}

\newtheorem{theorem}{Theorem}[section]
\newtheorem{lemma}[theorem]{Lemma}
\newtheorem{proposition}[theorem]{Proposition}
\newtheorem{corollary}[theorem]{Corollary}

\theoremstyle{definition}
\newtheorem{definition}[theorem]{Definition}
\newtheorem{example}[theorem]{Example}

\theoremstyle{remark}
\newtheorem{remark}[theorem]{Remark}
\numberwithin{equation}{section}

\title{Hyperbolicity in unbounded convex domains}
\author[F. Bracci]{Filippo Bracci}
\address{F. Bracci: Dipartimento Di Matematica\\
Universit\`{a} di Roma \textquotedblleft Tor Vergata\textquotedblright\ \\
Via Della Ricerca Scientifica 1, 00133 \\
Roma, Italy} \email{fbracci@mat.uniroma2.it}

\author[A. Saracco]{Alberto Saracco}
\address{A. Saracco: Scuola Normale Superiore \\
Piazza dei Cavalieri 7, 56126 \\
Pisa, Italy} \email{a.saracco@sns.it}
\thanks{The second author was partially supported by PRIN project \lq\lq Proprietà geometriche delle varietà reali e complesse\rq\rq}




\begin{abstract} We provide several equivalent
characterizations
of Kobayashi hyperbolicity in unbounded convex domains in terms
of peak and anti-peak functions at infinity, affine lines,
Bergman metric and iteration theory.
\end{abstract}

\maketitle

\section{Introduction}

Despite the fact that linear convexity is not an invariant
property in complex analysis, bounded convex domains in $\C^N$
have been very much studied as  prototypes for  the general
situation.

In particular, by Harris' theorem \cite{Ha} (see also,
\cite{Abate}, \cite{Kob}) it is known that bounded convex
domains are always Kobayashi complete hyperbolic (and thus by
Royden's theorem, they are also taut and hyperbolic). Moreover,
by Lempert's theorem \cite{Le}, \cite{Le1}, the Kobayashi
distance can be realized by means of extremal discs. These are
the basic cornerstone for many useful results, especially in
pluripotential theory and iteration theory.

On the other hand, not much is known about {\sl unbounded}
domains. Clearly, the geometry at infinity must play some
important role. In this direction, Gaussier \cite{Gau} gave
some conditions in terms of existence of {\sl peak and
anti-peak} functions at infinity for an unbounded domain to be
hyperbolic, taut or complete hyperbolic. Recently, Nikolov and
Pflug \cite{NP} deeply studied conditions at infinity which
guarantee hyperbolicity, up to a characterization of
hyperbolicity in terms of the asymptotic behavior of the
Lempert function.

In these notes we restrict ourselves to the case of {\sl
unbounded convex} domains, where, strange enough, many open
questions in the previous directions seem to be still open. In
particular an unbounded convex domain needs not to be
hyperbolic, as the example of $\C^k$ shows. Some estimates on the Caratheodory and Bergman metrics in convex domains were obtained by Nikolov and Pflug in \cite{NP03a}, \cite{NP03b}. The question is whether one can understand easily hyperbolicity of
unbounded convex domains in terms of geometric or analytic
properties. A result in this direction was obtained by Barth \cite{Barth}, who proved the equivalence of properties (1), (2) and (6) in the theorem below.

The aim of the present paper is to show that actually for
unbounded convex domains, hyperbolicity can be characterized in
many different ways and can be  easily inferred just looking at
a single boundary point.

The dichotomy we discovered for unbounded convex domains is
rather stringent: either the domain behaves like a bounded
convex domain or it behaves like $\C^k$. In particular, this
provides examples of unbounded domains which admit the Bergman
metric and are complete with respect to it.

The main result of these notes is the following (notations and
terminology  are standard and will be recalled in the next
section):

\begin{theorem}\label{main}
Let $D\subset \C^N$ be a (possibly unbounded) convex domain.
The following are equivalent:
\begin{enumerate}
\item $D$ is biholomorphic to a bounded domain;
\item $D$ is (Kobayashi) hyperbolic; 	
\item $D$ is taut; 	
\item $D$ is complete (Kobayashi) hyperbolic; 	
\item $D$ does not contain nonconstant entire curves; 	
\item $D$ does not contain complex affine lines; 	
\item $D$ has $N$ linearly independent separating real
hyperplanes;
\item $D$ has peak and antipeak functions (in the sense of Gaussier) at
infinity;
\item $D$ admits the Bergman metric $b_D$.
\item $D$  is complete with respect to the  Bergman metric
$b_D$.
\item for any $f:D\to D$ holomorphic such that the sequence of its iterates $\{f^{\circ k}\}$ is not
compactly divergent there exists $z_0\in D$ such that
$f(z_0)=z_0$.
\end{enumerate}
\end{theorem}

The first implications of the theorem allow to obtain the
following {\sl canonical complete hyperbolic decomposition} for
unbounded convex domains, which is used in the final part of
the proof of the theorem itself.

\begin{proposition}\label{standard-decomp}
Let $D\subset \C^N$ be a (possibly unbounded) convex domain.
Then there exist a unique $k$ ($0\leq k\leq N$) and
 a unique  complete hyperbolic convex domain $D'\subset\C^k$, such that,
up to a linear change of coordinates, $D\ =\ D'\times\C^{N-k}.$
\end{proposition}

By using such a canonical complete hyperbolic decomposition,
one sees for instance that the ``geometry at infinity'' of an
unbounded convex domain can be inferred from the geometry of
any  finite point of its boundary (see the last section for
precise statements). For example, as an application of
Corollary \ref{mazzao} and Theorem \ref{main}, existence of
peak and anti-peak functions (in the sense of Gaussier) for an
unbounded convex domain equals  the absence of complex line in
the ``CR-part'' of the boundary of the domain itself. This
answers a question in Gaussier's paper (see \cite[pag.
 115]{Gau}) about geometric conditions for the existence in convex domains of
 peak and anti-peak plurisubharmonic functions at infinity.

The authors want to sincerely thank prof. Nikolov for helpful
conversations, and in particular for sharing his idea of
constructing antipeak functions.

\section{Preliminary}

A {\sl convex domain} $D\subset \C^N$ is a domain such that for
any couple $z_0, z_1\in D$ the real segment joining $z_0$ and
$z_1$ is contained in $D$. It is well known that for any point
$p\in\de D$ there exists (at least) one {\sl real separating
hyperplane} $H_p=\{z\in\C^n : \Re L(z)=a\}$, with $L$ a complex
linear functional and $a\in \R$ such that $p\in H_p$ and $D\cap
H_p=\emptyset$. Such a hyperplane $H_p$ is sometimes also
called a {\sl tangent hyperplane} to $D$ at $p$. We say that
$k$ separating hyperplanes $H_j=\{\Re L_j(z)=a_j\}$,
$j=1,\ldots, k$, are linearly independent if $L_1,\ldots, L_k$
are linearly independent linear functionals.

Let $\D:=\{\zeta\in\C: |\zeta|<1\}$ be the unit disc. Let
$D\subset \C^N$ be a domain. The Kobayashi pseudo-metric for
the point $z\in D$ and vector $v\in\C^N$ is defined as
\[
\kappa_D(z;v):=\inf\{\lambda>0| \exists
\varphi:\D\stackrel{\hbox{hol}}{\longrightarrow} D ,
\varphi(0)=z, \varphi'(0)=v/\lambda\}.
\]
If $\kappa_D(z;v)>0$ for all $v\neq 0$ then $D$ is said to be
{\sl (Kobayashi) hyperbolic}. The pseudo-distance $k_D$
obtained by integrating $\kappa_D$ is called the {\sl Kobayashi
pseudodistance}. The domain $D$ is {\sl (Kobayashi) complete
hyperbolic} if $k_D$ is complete.

The Carath\'eodory pseudo-distance $c_D$ is defined by
\[
c_D(z,w)=\sup\{k_{\D}(f(z),f(w)): f:D\to \D \
\hbox{holomorphic}\}.
\]
In general, $c_D\leq k_D$.

We refer the reader to the book of Kobayashi \cite{Kob} for
properties of Kobayashi and Carath\'eodory metrics and
distances.

Another (pseudo)distance that can be introduced on the domain
$D$ is the {\sl Bergman (pseudo) distance} (see, {\sl e.g.},
\cite[Sect. 10, Ch. 4]{Kob}). Let $\{e_j\}$ be a orthonormal
complete basis of the space of square-integrable holomorphic
functions on $D$. Then let
\[
l_D(z,\overline{w}):=\sum_{j=0}^\infty e_j(z)
\overline{e}_j(\overline{w}).
\]
If $l_D(z,\overline{z})>0$ one can define a symmetric form
$b_D:=2\sum h_{jk} dz_j\otimes d\overline{z}_k$, with
$h_{jk}=\frac{\de^2\log b_D(z,\overline{z})}{\de z_j, \de
\overline{z}_k}$, which is a positive semi-definite Hermitian
form, called the {\sl Bergman pseudo-metric} of $D$. If $b_D$
is positive definite everywhere, one says that $D$ {\sl admits
the Bergman metric $b_D$}. For instance, $\C^k$, $k\geq 1$ does
not support square-integrable holomorphic functions, therefore
$l_{\C^k}\equiv 0$ and  $\C^k$ does not admit the Bergman
metric.

For the next result, see \cite[Corollaries 4.10.19,
4.10.20]{Kob}:
\begin{proposition}\label{bergman}
Let $D\subset \C^N$ be a domain.
\begin{enumerate}
\item Assume $b_D(z,\overline{z})>0$ for all $z\in D$. If $c_D$ is a distance, if it induces
the topology of $D$ and if the $c_D$-balls are compact, then
$D$ admits the Bergman metric $b_D$ and it is complete with
respect to $b_D$.
\item If  $D$ is a  bounded convex domain then it (admits
the Bergman metric and it) is complete with respect to the
Bergman metric.
\end{enumerate}
\end{proposition}

Let $G$ be another domain. We recall that if $\{\v_k\}$ is a
sequence of holomorphic mappings from $G$ to $D$, then the
sequence is said to be {\sl compactly divergent} if for any two
compact sets $K_1\subset G$ and $K_2\subset D$ it follows that
$\sharp\{k\in \N: \v_k(K_1)\cap K_2\neq \emptyset\}<+\infty$.

A family $\mathcal F$ of holomorphic mappings from $G$ to $D$
is said to be {\sl normal} if each sequence of $\mathcal F$
admits a subsequence which is either compactly divergent or
uniformly convergent on compacta.

If the family of all holomorphic mappings from the unit disc
$\D$ to $D$ is normal, then $D$ is said to be {\sl taut}. It is
known:

\begin{theorem}
Let $D\subset \C^N$ be a domain.
\begin{enumerate}
\item(Royden) $D\ \hbox{complete hyperbolic} \Rightarrow D\ \hbox{taut}
\Rightarrow D\ hyperbolic$.
\item(Kiernan) If $D$ is bounded then $D$ is hyperbolic.
\item(Harris) If $D$ is a bounded convex domain then $D$ is
complete hyperbolic.
\end{enumerate}
\end{theorem}

The notion of (complete) hyperbolicity is pretty much related
to existence of peak functions at each boundary point. In case
$D$ is an unbounded domain, H. Gaussier \cite{Gau} introduced
the following concepts of ``peak and antipeak functions'' at
infinity,  which we use in the sequel:

\begin{definition}
A function $\varphi:\overline D\to\R\cup\{-\infty\}$ is called
a \emph{global peak plurisubharmonic function at infinity} if
it is plurisubharmonic on $D$, continuous up to $\overline D$
(closure in $\C^N$) and
$$\begin{cases}\lim\limits_{z\to\infty}\varphi(z)=0, \\ \varphi(z)<0 \quad \forall z\in \overline{D}.\end{cases}$$
A function $\varphi:\overline D\to\R\cup\{-\infty\}$ is called
a \emph{global antipeak plurisubharmonic  function at infinity}
if it is plurisubharmonic on $D$, continuous up to $\overline
D$ and
$$\begin{cases}\lim\limits_{z\to\infty}\varphi(z)=-\infty,
\\ \varphi(z)>-\infty \quad \forall z\in \overline{D}.\end{cases}$$

For short we will simply call them {\sl peak and antipeak
functions (in the sense of Gaussier) at infinity}.
\end{definition}

Gaussier proved the following result:

\begin{theorem}[Gaussier]
Let $D\subset \C^N$ be an unbounded domain. Assume that $D$ is
locally taut at each point of $\de D$ and there exist peak and
antipeak functions (in the sense of Gaussier) at infinity. Then
$D$ is taut.
\end{theorem}

Obviously a convex domain is locally taut at each boundary
point, thus  tautness  follows from existence of peak and
antipeak functions (in the sense of Gaussier) at infinity.

Finally, if $f:D\to D$ is a holomorphic function, the sequence
of its iterates $\{f^{\circ k}\}$ is defined by induction as
$f^{\circ k}:=f^{\circ (k-1)}\circ f$. If $f$ has a fixed point
$z_0\in D$, then $\{f^{\circ k}\}$ is not compactly divergent.
On the other hand, depending on the geometry of $D$, there
exist  examples of holomorphic maps $f$ such that $\{f^{\circ
k}\}$ is not compactly divergent but $f$ has no fixed points in
$D$. It is known (see \cite{Abate2}) that

\begin{theorem}[Abate]\label{abaabaye}
Let $D\subset \C^N$ be a taut domain. Assume that
$H^j(D;\mathbb Q)=0$ for all $j>0$ and let $f:X\to X$
holomorphic. Then the sequence of iterates $\{f^{\circ k}\}$ is
compactly divergent if and only if $f$ has no periodic points
in $D$.

If $D$ is a bounded convex domain then the sequence of iterates
$\{f^{\circ k}\}$ is compactly divergent if and only if $f$ has
no fixed points in $D$.
\end{theorem}

\section{The proof of Theorem \ref{main}}

The proof of Theorem \ref{main} is obtained in several steps,
which might be of some interest by their own.

For a domain $D\subset\C^N$ let us denote by $\delta_D$ the
{\sl Lempert function} given by
\[
\delta_D(z,w)=\inf\{\omega(0,t): t\in (0,1), \exists \v\in \Hol
: \v(0)=z, \v(t)=w\}.
\]
The Lempert function is not a pseudodistance in general because
it does not enjoy the triangle inequality. The Kobayashi
pseudodistance is the largest minorant of $\delta_D$ which
satisfies the triangle inequality. The following lemma (known
as Lempert's theorem in case of bounded convex domains) is
probably known, but we provide its simple proof due to the lack
of reference.

\begin{lemma}\label{Lempert}
Let $D\subset \C^N$ be a (possibly unbounded) convex domain.
Then $k_D=\delta_D=c_D$.
\end{lemma}

\begin{proof}
The result is due to Lempert \cite{Le1} in case $D$ is bounded.
Assume $D$ is unbounded. Let $D_R$ be the intersection of $D$
with a ball of center the origin and radius $R>0$. For $R>>1$
the set $D_R$ is a nonempty convex bounded domain. Therefore
$k_{D_R}=\delta_{D_R}=c_{D_R}$. Now $\{D_R\}$ is an increasing
sequence of domains whose union is $D$. Hence,  $\lim_{R\to
\infty} k_{D_R}=k_D$, $\lim_{R\to \infty} c_{D_R}=c_D$ and
$\lim_{R\to \infty} \delta_{D_R}=\delta_D$ (see, {\sl e.g.},
\cite[Prop. 2.5.1]{JP} and \cite[Prop. 3.3.5]{JP}). Thus
$k_D=\delta_D=c_D$.
\end{proof}

\begin{proposition}\label{palleconvesse} Let $D\subset\C^N$ be a (possibly unbounded) convex domain. Then the Kobayashi balls in
$D$ are convex.
\end{proposition}
\begin{proof} For the bounded case, see \cite[Proposition 2.3.46]{Abate}. For the unbounded case,
let $B_\epsilon$ be the Kobayashi ball of radius $\epsilon$ and
center $z_0\in D$, let $D_R$ be the intersection of $D$ with an
Euclidean ball of center the origin and  radius $R>0$, and let
$B^R_\epsilon$ be the Kobayashi ball  in $D_R$  of radius
$\epsilon$ and center $z_0$. Then the convex sets
$B^R_\epsilon\subset B^{R+\delta}_\epsilon\subset B_\epsilon$
for all $R>>1$, $\delta>0$, and their convex increasing union
$\cup_R B^R_\epsilon=B_\epsilon,$ since $\lim_{R\to \infty}
k_{D_R}=k_D$.
\end{proof}

\begin{lemma}\label{extremal}
Let $D\subset \C^N$ be a (possibly unbounded) taut convex
domain. Then for any couple $z,w\in D$ there exists $\v \in
\Hol$ such that $\v(0)=z, \v(t)=w$ $t\in [0,1)$ and
$k_D(z,w)=\omega(0,t)$.
\end{lemma}

\begin{proof}
By Lemma \ref{Lempert}, $k_D=\delta_D$, so there exists a
sequence $\{\v_k\}$ of holomorphic discs and $t_k\in (0,1)$
such that $\v_k(0)=z$ and $\v_k(t_k)=w$ and
\[
k_D(z,w)=\lim_{k\to\infty}\omega(0,t_k).
\]
Since $D$ is taut and $\v_k(0)=z$ for all $k$, we can assume
that $\{\v_k\}$ converges uniformly on compacta to a
(holomorphic) map $\v:\D\to D$. Then $\v(0)=z$. Moreover, since
$k_D(z,w)<\infty$,   there exists $t_0<1$ such that $t_k\leq
t_0$ for all $k$. We can assume (up to subsequences) that
$t_k\to t\leq t_0$. Then
\[
k_D(z,w)=\lim_{k\to\infty}\omega(0,t_k)=\omega(0,t).
\]
Moreover, $\v(t)=\lim_{k\to\infty}\v_k(t_k)=w$ and we are done.
\end{proof}

\begin{proposition}\label{TiifCH}
Let $D\subset \C^N$ be a (possibly unbounded) convex domain.
Then $D$ is taut if and only if it is complete hyperbolic.
\end{proposition}

\begin{proof}
One direction is contained in Royden's theorem. Conversely,
assume that $D$ is taut. We are going to prove that every
closed Kobayashi balls is compact (which is equivalent to be
complete hyperbolic, see \cite{Kob} or \cite[Proposition
2.3.17]{Abate}).

Let $R>0$, $z\in D$ and let $B(z,R)=\{w\in D: k_D(z,w)\leq
R\}$. If $B(z,R)$ is not compact then there exists a sequence
$\{w_k\}$ such that $w_k\to p\in \de D\cup\{\infty\}$ and
$k_D(z,w_k)\leq R$. For any $k$, let $\v_k\in \Hol$ be the
extremal disc given by Lemma \ref{extremal} such that
$\v_k(0)=z, \v_k(t_k)=w_k$ for some $t_k\in (0,1)$ and
$k_D(z,w_k)=\omega(0,t_k)$.

Notice that, since $k_D(z,w_k)\leq R$, then there exists $t_0<1$
such that $t_k\leq t_0$ for all $k$. We can assume up to
subsequences that $t_k\to t$ with $t<1$. Since $D$ is taut and
$\v_k(0)=z$, up to extracting subsequences, the sequence
$\{\v_k\}$ is converging uniformly on compacta to a holomorphic
disc $\v:\D\to D$ such that $\v(0)=z$. However,
\[
\v(t)=\lim_{k\to \infty}\v_k(t_k)=\lim_{k\to\infty} w_k=p,
\]
a contradiction. Therefore $B(z,R)$ is compact and $D$ is
complete hyperbolic.
\end{proof}

For the next proposition, cfr. \cite[Lemma 3]{Dr}.

\begin{proposition}\label{Siegel} Let $D\subset \C^N$ be a convex domain, which does not contain complex affine
lines. Then there exist $\{L_1=0\},\ldots,\{L_N=0\}$ linearly
independent  hyperplanes containing the origin and $a_1,\ldots,
a_N\in\R$ such that
$$D \ \subset\ \{\Re\,L_1>a_1,\ldots,\Re\,L_N>a_N\}.$$
\end{proposition}
\begin{proof} Without loss of generality we can assume that $O\in D$. Since $D$ does not contain
complex affine lines, $\de D$ is not empty. Take a point
$p_1\in \de D$ and a tangent real hyperplane through $p_1$,
given by $\{\Re\,L_1=a_1\}$ (if the boundary is smooth there is
only one tangent hyperplane), where $L_1$ is defined so that
$D\subset\{\Re\,L_1>a_1\}$.

Suppose that $L_1,\ldots,L_k$, $k<N$, are already defined, they are linearly independent and
$$D \ \subset\ \{\Re\,L_1>a_1,\ldots,\Re\,L_k>a_k\}.$$
 The  intersection
$h_k=\cap_1^k \{L_i=0\}$ is a complex $(N-k)$-dimensional plane
through the origin $O$ (which is also contained in $D$ by
hypothesis). Since $D$ does not contain complex affine lines,
$\de D\cap h_k$ is not empty. Take a point $p_{k+1}\in \de
D\cap h_k$ and consider a tangent real hyperplane through
$p_{k+1}$, $\{\Re\,L_{k+1}=a_{k+1}\}$, where $L_{k+1}$ is
defined so that $D\subset\{\Re\,L_{k+1}>a_{k+1}\}$. By
construction $L_{k+1}$ is linearly independent from
$L_1,\ldots,L_k$ and
$$D \ \subset\ \{\Re\,L_1>a_1,\ldots,\Re\,L_{k+1}>a_{k+1}\}.$$
Continuing this way, the proof is  concluded.
\end{proof}

Now we are in a good shape to prove part of Theorem \ref{main}:

\begin{proposition}\label{equiv} Let $D\subset \C^N$ be a convex domain. The following are equivalent:
\begin{enumerate}
\item $D$ is biholomorphic to a bounded domain;
\item $D$ is (Kobayashi) hyperbolic; 	
\item $D$ is taut; 	
\item $D$ is complete (Kobayashi) hyperbolic; 	
\item $D$ does not contain nonconstant entire curves; 	
\item $D$ does not contain complex affine lines; 	
\item $D$ has $N$ linearly independent separating real
hyperplanes;
\item $D$ has peak and antipeak functions (in the sense of Gaussier) at
infinity;
\end{enumerate}
\end{proposition}

\begin{proof} (1) $\Rightarrow$ (2): every bounded domain in $\C^N$ is
hyperbolic by \cite{Ki} (see, also, \cite[Thm. 2.3.14]{Abate})

(2) $\Rightarrow$ (5) $\Rightarrow$ (6): obvious.

(6) $\Rightarrow$ (7): it is Proposition \ref{Siegel}.

(7) $\Rightarrow$ (1): let $L_1,\ldots, L_N$ be linearly
independent complex linear functionals and let $a_1,\ldots,
a_N\in \R$ be such that $\{\Re L_j=a_j\}$ for $j=1,\ldots, N$
are real separating hyperplanes for $D$. Up to sign changes, we
can assume that $D\subset \{\Re L_j>a_j\}$. Then the map
\[
F(z_1,\ldots, z_N):=\left(\frac{1}{L_1(z)-a_1+1}, \ldots,
\frac{1}{L_N(z)-a_N+1}\right)
\]
maps $D$ biholomorphically on a bounded convex domain of
$\C^N$.

(6) $\Rightarrow$ (8): let $L_1,\ldots,L_N$ be as in Proposition \ref{Siegel}. Up to a linear change of coordinates, we can suppose that $z_j=L_j$ for all $1\leq j\leq N$. A peak function is given by
$$-\Re\,\sum_{j=1}^N \frac1{z_j-a_j+1}.$$
Let $D_j:= \{\Re L_j>a_j\}$. Then $D\subset\prod_{j=1}^N D_j$,
and $D_j$ is biholomorphic to $\D$ for each $j$. In particular
$\C\setminus D_j$ is not a polar set. We may assume that
$0\not\in D_j.$ Let $G_j$ be the image of $D_j$ under the
transformation $z\to1/z.$ Since $\C\setminus G_j$ is not a
polar set, there exists $\varepsilon>0$ such that $\C\setminus
G_j^\varepsilon$ is not polar, too, where
$G_j^\varepsilon=G_j\cup\varepsilon\D.$ Denote by
$g_j^\varepsilon$ the Green function of $G_j^\varepsilon.$ Then
$h_j=g_j^\varepsilon(0;\cdot)$ is a negative harmonic function
on $G_j$ with $\lim_{z\to 0}h_j(z)=-\infty$ and
$\inf_{G_j\setminus r\D}h_j>-\infty$ for any $r>0.$ Then
$\psi_j(z)=h_j(1/z)$ is an antipeak function of $D_j$ at
$\infty$ and hence $\psi=\sum_{j=1}^{N'}\psi_j$ is an antipeak
function for $D$ at $\infty.$

(8) $\Rightarrow$ (3): it is  Gaussier's theorem \cite[Prop.
2]{Gau}.

(3) $\Rightarrow$ (4): it is Proposition
\ref{TiifCH}.

(4) $\Rightarrow$ (3) $\Rightarrow$ (2): it is  Royden's
theorem \cite[Prop. 5, pag. 135 and Corollary p.136]{Ro}.
\end{proof}

As a consequence we have Proposition \ref{standard-decomp},
which gives  a \emph{canonical complete hyperbolic
decomposition} of a convex domain as the product of a complete
hyperbolic domain and a copy of $\C^k$.

\begin{proof}[Proof of Proposition \ref{standard-decomp}] We prove the result by induction on $N$. If $N=1$ then either $D=\C$ or $D$ is biholomorphic to the disc
and hence (complete) hyperbolic.

Assume the result is true for $N$, we prove it holds for $N+1$.
Let $D\subset \C^{N+1}$ be a convex domain. Then, by
Proposition \ref{equiv}, either $D$ is complete hyperbolic or
$D$ contains an affine line, say, up to a linear change of
coordinates
$$l_{N+1}=\{z_1,\ldots, z_N=0\}\subset D.$$
Clearly, there exists $c\in \C$ such that
$D\cap\{z_{N+1}=c\}\neq \emptyset$. Up to translation we can
assume $c=0$. Let us define
$$D_N=D\cap\{z_{N+1}=0\}.$$
$D_N\subset\C^N$ is convex. We claim that $D=D_N\times\C$.
Induction will then conclude the proof.

Let $z_0\in D_N$. We want to show that $(z_0,\zeta)\in D$ for
all $\zeta\in\C$. Since $l_{N+1}\subset D$ then $(0,\zeta)\in
D$ for all $\zeta\in \C$. Assume $z_0\neq 0$. Fix $\zeta\in\C$.
Since $D_N$ is open, there exists $\varepsilon_0>0$ such that
$z_1:=(1+\varepsilon_0) z_0\in D_N$. Since $D$ is convex, for
any $t\in [0,1]$ it follows $t(z_1,0)+(1-t)(0,\xi)\in D$ for
all $\xi\in\C$. Setting
$\xi_0:=\frac{1+\varepsilon_0}{\varepsilon_0}\zeta\in\C$ and
$t_0=(1+\varepsilon_0)^{-1}\in(0,1)$ we obtain
$$(z_0,\zeta)=t_0\left(z_1,0\right)\
+\ (1-t_0)\left(0,\xi_0\right)\in D,$$ completing the proof.
\end{proof}

In order to finish the proof of Theorem \ref{main} we need to
show that the first eight conditions, which are all and the
same thanks to Proposition \ref{equiv}, are equivalent to (9),
(10), (11).

\begin{proof}[Proof of Theorem \ref{main}] Conditions (1) to
(8)  are all equivalent by Proposition \ref{equiv}.

(10)$\Rightarrow$(9): obvious.

(4)$\Rightarrow$(10): By Lemma \ref{Lempert}, the Caratheodory
distance $c_D$ equals the Kobayashi distance $k_D$, thus, since
$D$ is (Kobayashi) complete hyperbolic,  $c_D$ is a distance
which induces the topology on $D$ and the $c_D$-balls are
compact. By Proposition \ref{bergman} then $D$ admits the
Bergman metric and it is complete with respect to it.

(9)$\Rightarrow$(4): Assume $D$ is not complete hyperbolic.
Then by Proposition \ref{standard-decomp}, up to a linear
change of coordinates, $D=D'\times \C^k$ for some complete
hyperbolic domain $D'$ and $k\geq 1$. By the product formula
(see \cite[Prop. 4.10.17]{Kob}) $l_D=l_{D'}\cdot l_{\C^k}\equiv
0$ and thus $D$ does not admit the Bergman metric.

(11) $\Rightarrow$ (4): Assume $D$ is not complete hyperbolic.
We have to exhibit a holomorphic self-map $f:D\to D$ such that
$\{f^{\circ k}\}$ is not compactly divergent but there exists
no $z_0\in D$ such that $f(z_0)=z_0$.

By Proposition \ref{standard-decomp}, up to a linear change of
coordinates, $D=D'\times \C^k$ for some complete hyperbolic
domain $D'$ and $k\geq 1$. Let
\[
f: D'\times \C^{k-1}\times \C \ni (z,w',w)\mapsto (z,w',
e^w+w)\in D'\times \C^{k-1}\times \C.
\]
Then clearly $f$ has no fixed points in $D$. However, if
$w_0=\log(i\pi)$, then $f^{\circ 2}(z,w',w_0)=(z,w',w_0)$, and
therefore the sequence $\{f^{\circ k}\}$ is not compactly
divergent.

(4) $\Rightarrow$ (11): According to the theory developed so
far, if $D$ is complete hyperbolic, then it is taut and its
Kobayashi balls are convex and compact. With these ingredients,
the proof for {\sl bounded} convex domains go through also in
the unbounded case (see \cite[Thm. 2.4.20]{Abate}).
\end{proof}

\section{Applications}

\begin{corollary}
Let $D\subset \C^N$ be a convex domain. If there exists a point
$p\in\de D$ such that $\de D$ is strongly convex at $p$ then
$D$ is complete hyperbolic.
\end{corollary}

\begin{proof}
By Proposition \ref{standard-decomp}, if $D$ were not complete
hyperbolic, up to linear changes of coordinates, $D=D'\times
\C^k$ for some complete hyperbolic convex domain $D'$ and
$k\geq 1$. Then $\de D=\de D'\times \C^k$ could not be strongly
convex anywhere.
\end{proof}

Note that the converse to the previous corollary is false: the
half-plane $\{\zeta\in \C: \Re \zeta>0\}$ is a complete
hyperbolic convex domain in $\C$ with boundary which is nowhere
strongly convex.

\begin{proposition}
Let $D\subset\C^N$ be an (unbounded) convex domain. Then $D$
has canonical complete hyperbolic decomposition (up to a linear
change of coordinates) $D=D'\times\C^k$ with $D'$ complete
hyperbolic, if and only if for every $p\in
\partial D$ and every separating hyperplane $H_p$,  $(H_p\cap\partial D)\cap
i(H_p\cap\partial D)$ contains a copy of $\C^k$ but contains no
copies of $\C^{k+1}$.
\end{proposition}
\begin{proof} ($\Rightarrow$) Since $D=D'\times\C^k$, for every $p\in \partial D$ and every separating
hyperplane $H_p$,
$$(H_p\cap\partial D)\cap i(H_p\cap\partial D)=[(H_p\cap\partial D')\cap i(H_p\cap\partial D')]\times \C^k.$$
Since $D'$ is complete hyperbolic,  its boundary does not
contain complex lines.

($\Leftarrow$) Since $D$ is convex, $D=D'\times\C^{k'}$, by
Proposition \ref{standard-decomp}. By the first part of the
present proof, for every $p\in
\partial D$ and every separating hyperplane $H_p$,
$$(H_p\cap\partial D)\cap i(H_p\cap\partial D)= \C^{k'}.$$
Hence $k'=k$.
\end{proof}

\begin{corollary}\label{mazzao} Let $D\subset\C^N$ be an (unbounded) convex domain. If there exist  $p\in
\partial D$ and a separating hyperplane $H_p$ such that
$(H_p\cap\partial D)\cap i(H_p\cap\partial D)$ does not contain
any complex affine line then $D$ is complete hyperbolic.
Conversely, if $D$ is complete hyperbolic, then for any point
$p\in
\partial D$ and any separating hyperplane $H_p$, it follows
that $(H_p\cap\partial D)\cap i(H_p\cap\partial D)$ does not
contain any complex affine line.
\end{corollary}

As a final remark, we notice that, as Abate's theorem
\ref{abaabaye} is the cornerstone to the study of iteration
theory in bounded convex domains, our Theorem \ref{main} and
Proposition \ref{standard-decomp} can be used effectively well
to the same aim for unbounded convex domains. In fact, if
$D=D'\times \C^k$ is the canonical complete hyperbolic
decomposition of $D$, then a holomorphic self map $f:D\to D$
can be written in the coordinates $(z,w)\in D'\times\C^k$ as
$f(z,w)=(\v(z,w),\psi(z,w))$, where $\v: D'\times\C^k\to D'$
and $\psi: D'\times\C^k\to \C^k$. In particular, since $D'$ is
complete hyperbolic, then $\v$ depends only on $z$, namely,
$f(z,w)=(\v(z), \psi(z,w))$. The map $\psi(z,w)$ can be as
worse as entire functions in $\C^k$ are, but the map $\v$ is a
holomorphic self-map of a complete hyperbolic convex domains
and its dynamics goes   similarly to  that of holomorphic
self-maps of bounded convex domains. For instance, if the
sequence $\{f^{\circ k}\}$ is non-compactly divergent, then $f$
might have  no   fixed points, but the sequence $\{\v^{\circ
k}\}$ must have at least one.

\end{document}